\documentclass[11pt]{article}
\usepackage{a4wide}

\usepackage[numbers]{natbib}  % use [authoryear] for author-date style

\usepackage[fleqn]{amsmath}
\usepackage{amssymb}
\usepackage{tabularx}
\usepackage{ltablex}
\usepackage{csquotes}
\usepackage{graphicx}
\usepackage{placeins}
\usepackage{eurosym}
\keepXColumns
\usepackage[dvipsnames]{xcolor}
\usepackage{array} % Needed for new column type

\newcolumntype{C}[1]{>{\centering\arraybackslash}m{#1}}

\title{Discounting Approaches in Multi-Year Investment Modelling for Energy Systems}

\author{
Ni Wang\textsuperscript{1}, Diego A. Tejada-Arango\textsuperscript{1,2} \\
\parbox{\textwidth}{
\centering
\textit{\small
\textsuperscript{1}Energy \& Materials Transition, TNO, The Netherlands \\
\textsuperscript{2}Instituto de Investigación Tecnológica, Escuela Técnica Superior de Ingeniería,\\
Universidad Pontificia Comillas, Spain}
}
}
\date{}

\begin{document}

\maketitle

\begin{abstract}
This paper reviews discounting approaches for modeling multi-year energy investments, focusing on total versus annualised cost formulations. We discuss how time value of money is handled, and how salvage value and milestone-year weighting can address mismatches between asset lifetimes and model horizons. These methods are implemented in the open-source TulipaEnergyModel to support transparent and tractable long-term energy system planning.
\end{abstract}

\section{Introduction} \label{sec:intro}
This document is an updated and extended version of \cite{Tejada-Arango2023c},\footnote{This version supersedes the earlier report by Tejada-Arango et al. \cite{Tejada-Arango2023c}, providing a more detailed treatment of discounting approaches and a clearer narrative structure.} featuring a refined narrative structure and expanded treatment of discounting approaches.  It provides a comprehensive overview of the main methodologies for modeling multi-year investments in energy planning models, with a particular focus on how cost discounting is handled. To support a clear understanding of the underlying modeling logic, the document begins with conceptual explanations and gradually introduces mathematical formulations at varying levels of detail. Illustrative examples are included throughout to enhance clarity and highlight key modeling trade-offs.

\subsection{Common ground}

Before presenting the basic formulations, it is essential to establish some conceptual foundations. First, the term “multi-year investments” refers to investment decisions made at different points in time, typically measured in years. This reflects real-world investment behavior—for example, wind capacity may be added in year 0 and additional capacity in year 5. In this document, we simplify by assuming no lead or construction time; that is, the year of the investment decision is also the year in which the asset becomes operational.

Second, we consider the time value of money, which recognizes that the value of money depends on when it is received or spent. The standard approach is to express all future costs and revenues in terms of present value, typically referenced to year 0. This requires discounting future values to account for the opportunity cost of capital.

Two discounting concepts are commonly used: the Weighted Average Cost of Capital (WACC) and the social discount rate $R$. WACC is technology-specific, incorporating the financial risk and cost of capital associated with a particular asset. In contrast, $R$ reflects the broader societal preference for present over future consumption. While the two rates can differ, WACC generally incorporates the social discount rate as one of its components; hence, it is commonly assumed that $WACC \ge R$ . In this document, we adopt this assumption: $R$ is applied when discounting general future values, whereas WACC is used when evaluating technology-specific investment streams.

These two concepts—multi-year investment timing and discounting—are central to the discussions that follow. For clarity, we will revisit them throughout the document wherever misinterpretation may arise.

\subsection{Basic concepts of two approaches}

In the literature, there are two approaches that represent the investment costs:

\begin{itemize}
    \item \textit{Total Cost}: Considering the total cost across the lifetime of the asset $C_{y}^{T}$ (i.e., overnight cost) for the investment year $y$. 
    \item \textit{annualised Cost}: Converting the total cost into annualised costs $C_{y}^{A}$. Since there is time value of the money, $C_{y}^{A}$ are valued differently now and in the future. And since we always discuss values now, $C_{y}^{A}$ will be discounted back to now.

\end{itemize}

To illustrate the concepts, let us consider the following data as an example for the investments made in the first year ($y=0$):

\begin{itemize}
    \item Total cost $C_{0}^{T}=100$ [\euro/MW]
    \item $WACC_{0}=2$ \%
    \item Economic lifetime $LT=5$ years
\end{itemize}

The total cost $C_{0}^{T}$ is 100.

The annualised costs $C_{0}^{A}$ is given in the equation below describing its relationship with $C_{0}^{T}$, following which the solution for $C_{0}^{A}$ is $C_{0}^{I}=20.80$.

\begin{equation}
\label{eq:example_annualised_cost}
    C_{0}^{T} = 100 = C_{0}^{A} \cdot \left(1+\frac{1}{(1+2 \%)^{1}}+\frac{1}{(1+2 \%)^{2}}+\frac{1}{(1+2 \%)^{3}}+\frac{1}{(1+2 \%)^{4}} \right)
\end{equation}

Equation \ref{eq:example_annualised_cost} can be generalized as

\begin{equation}
 \label{eq:total_cost_from_annualised_cost}
    C_{y}^{T} = \sum_{j = y}^{y+LT-1}\frac{1}{(1+WACC_{y})^{j - y}} \cdot C_{j}^{A} \quad \forall y
\end{equation}

One standard method is to assume that the annuities remain constant throughout the lifetime and the first year is not discounted, as explain in \cite{Centeno2018}. This means that we can replace $C_{j}^{A}$ with $C_{y}^{A}$. After some math operations, this formula becomes 

\begin{equation}
 \label{eq:annualised_cost_from_total_cost}
    C_{y}^{A} = \frac{WACC_{y}}{(1+WACC_{y}) \cdot \left( 1- \frac{1}{(1+WACC_{y})^{LT}} \right)} \cdot C_{y}^{T} \quad \forall y
\end{equation}

The reference \cite{Brown2020} provides a concise summary of the approach utilized in popular planning models within the energy sector, including TIMES, DIMENSION, SWITCH, and PRIMES.

\subsection{Basic formulations of the two approaches}

In this section, we revisit the concepts discussed above in the context of a multi-year investment planning model. Equation (\ref{eq:objective_function}) presents a general formulation of the objective function, while Equation (\ref{eq:production_limit}) captures how investment decisions constrain generation output over time, based on asset lifetimes. For clarity, the objective function omits several cost components such as fixed operation and maintenance costs, decommissioning costs, and mothballing costs, that are commonly included in more detailed formulations. Beyond these standard techno-economic factors, other relevant cost categories such as indirect costs \cite{Wang2023}, and market-related costs associated with evolving market designs and trading mechanisms \cite{Estanqueiro2022} can also be incorporated to enhance the realism and policy relevance of the model.

\begin{equation}
 \label{eq:objective_function}
    \min_{x_{y},p_{ykt}} \quad C^{I} + C^{O} 
\end{equation}

\begin{equation}
 \label{eq:production_limit}
 \begin{aligned}
    &s.t. \\
    &p_{ykt} \leq \sum_{j=\max(y-LT+1,0)}^{y}x_{j} \quad \forall y \forall k \forall t
 \end{aligned}
\end{equation}

Where $C^{I}$ refers to the investment cost and $C^{O}$ to the operational cost. A general definition for $C^{O}$ is shown in equation (\ref{eq:basic_operational_cost_formulation}):

\begin{equation}
 \label{eq:basic_operational_cost_formulation}
    C^{O} = \sum_{y=0}^{Y-1}\frac{1}{(1+R)^{y}} \cdot C_{y}^{op} \sum_{k}W_{yk}^{op} \cdot \sum_{t} p_{ykt} \\
\end{equation}

$C_{y}^{op}$ is the costs in the year $y$, which has to be discounted back to now only using R.

We now focus on the investment cost $C^{I}$.

\subsubsection{Total cost approach}
Equation (\ref{eq:basic_total_cost_formulation}) shows the total investment cost across the modelled years, considering the social discount rate $R$.

\begin{equation}
 \label{eq:basic_total_cost_formulation}
 \begin{aligned}
    & C^{I} = \sum_{y=0}^{Y-1}\frac{1}{(1+R)^{y}} \cdot C_{y}^{T}x_{y} \\
    & C^{1} = C_{0}^{T}x_{0} + 
    \frac{1}{(1+R)^{1}} \cdot C_{1}^{T}x_{1} + 
    \frac{1}{(1+R)^{2}} \cdot C_{2}^{T}x_{2} + \cdots
 \end{aligned}
\end{equation}

\subsubsection{Annualised cost approach}
Equation (\ref{eq:basic_annualised_cost_formulation}) is obtained by replacing (\ref{eq:total_cost_from_annualised_cost}) in (\ref{eq:basic_total_cost_formulation}). Bear in mind that the value of the parameter $C_{y}^{A}$ is determined by (\ref{eq:annualised_cost_from_total_cost}).

\begin{equation}
 \label{eq:basic_annualised_cost_formulation}
 \begin{aligned}
    C^{I} = & \sum_{y=0}^{Y-1} \frac{1}{(1+R)^{y}} \cdot C_{y}^{A} x_{y} \cdot \sum_{j = y}^{\min(y+LT-1,Y-1)}\frac{1}{(1+WACC_y)^{j-y}} \\
    C^{I} = & C_{0}^{A} \cdot x_{0} \cdot (1 + \frac{1}{(1+WACC_0)^{1}} + \frac{1}{(1+WACC_0)^{2}} + \cdots) + \\
            & \frac{1}{(1+R)^{1}} \cdot C_{1}^{A} \cdot x_{1} \cdot (1 + \frac{1}{(1+WACC_1)^{1}} + \frac{1}{(1+WACC_1)^{2}} + \cdots) + \\ 
            & \frac{1}{(1+R)^{2}} \cdot C_{2}^{A} \cdot x_{2} \cdot (1 + \frac{1}{(1+WACC_2)^{1}} + \frac{1}{(1+WACC_2)^{2}} + \cdots) + \\
            & \cdots
 \end{aligned}
\end{equation}

This equation can be explained as follows:
\begin{itemize}
    \item the annualised cost $C_{y}^{A}$ from the investment $x_y$ made in year $y$ is valued at year $y$.
    \item for the years after, its value has to be discounted from a future year to the year $y$ by WACC. It has to be emphasized here that we do not use R for this part of the equation because R only considers the social discount but not the technology-specific discount. The former applies when counting from year 0 to the year of the investment, i.e., year $y$. And the latter applies from the year of investment to the future, among which this investment is annualised.
    \item for the years before, since the investment has not yet been made, its value is discounted from the year $y$ to the year $0$ only using R. 
\end{itemize}

\section{Analysis}

Equations (\ref{eq:basic_total_cost_formulation}) and (\ref{eq:basic_annualised_cost_formulation}) calculate the same investment costs. Their equivalence can be checked by revisiting the definition of $C_y^T$ (Equation (\ref{eq:total_cost_from_annualised_cost})) and then comparing Equations (\ref{eq:basic_total_cost_formulation}) and (\ref{eq:basic_annualised_cost_formulation}) as follows.

\begin{equation*}
 \label{eq:equivalence_total_and_annualised_cost}
 \begin{aligned}
    C_{0}^{T} = & C_{0}^{A} \cdot \left(1 + \frac{1}{(1+WACC_{0})^{1}} + \frac{1}{(1+WACC_{0})^{2}} + \cdots \right) \\ 
    C_{1}^{T} = & C_{1}^{A} \cdot \left(1 + \frac{1}{(1+WACC_{1})^{1}} + \frac{1}{(1+WACC_{1})^{2}} + \cdots \right) \\ 
    C_{2}^{T} = & C_{2}^{A} \cdot \left(1 + \frac{1}{(1+WACC_{2})^{1}} + \frac{1}{(1+WACC_{2})^{2}} + \cdots \right) \\
    \vdots \qquad & \vdots
 \end{aligned}
\end{equation*}

Despite being mathematically equivalent, the total cost approach and the annualised approach are still different because of the starting points. In the total cost approach, one uses the parameter $C_y^T$, i.e., the overnight cost in the objective function. In the annualised cost approach, one uses the parameter $C_{y}^{A}$, i.e., the annualised cost in the objective function. We have shown that these data can be derived from one to the other using Equations (\ref{eq:total_cost_from_annualised_cost}) and (\ref{eq:annualised_cost_from_total_cost}) given the availability of all the used parameters. However, depending on the data availability, one may prefer one approach over the other.  

Apart from the data availability, the considered modelling horizon may also influence the choice of the approach. The above analysis has two assumptions. First, the modelling horizon is infinite, i.e., greater than the lifetime of the investments, such that all costs over the lifetime are considered. Second, investments are made yearly. However, the usual user case considers that modelled years do not cover the entire investment lifetime and use milestone years (or representative years) to reduce the computational burden. We will now analyse these cases and check how they differ from the ideal approaches.

\subsection{Salvage value}

Let us start with the case where modelled years do not cover the entire investment lifetime. We will show the formulations of the two approaches and explain a problem with the total cost approach. Then, to deal with the problem, we will introduce a concept: salvage value.  

Consider this example where the last modelled year ($Y$) is 2050, and we have an investment made in 2040 with a lifetime of 15 years. It follows that the last four years will not be considered in the optimisation. In this situation:

\begin{itemize}
    \item The annualised cost approach in (\ref{eq:basic_annualised_cost_formulation}) will only consider the annualised costs from 2040 to 2050. The ideal cost considering the lifetime (using equation (\ref{eq:basic_annualised_cost_formulation})) should be
    \begin{equation}
     \label{eq:example_annualised_cost_formulation_ideal}
     \begin{aligned}
                C^{I} = & C_{0}^{A} \cdot x_{0} \cdot (1 + \frac{1}{(1+WACC_0)^{1}} + \frac{1}{(1+WACC_0)^{2}} + \cdots + \frac{1}{(1+WACC_0)^{14}}) 
     \end{aligned}
    \end{equation}
    But because we only model till 2050, we use 
    \begin{equation}
     \label{eq:example_annualised_cost_formulation_real}
     \begin{aligned}
                C^{I} = & C_{0}^{A} \cdot x_{0} \cdot (1 + \frac{1}{(1+WACC_0)^{1}} + \frac{1}{(1+WACC_0)^{2}} + \cdots + \frac{1}{(1+WACC_0)^{10}}) 
     \end{aligned}
    \end{equation}
    
    Therefore, the annualised costs of the non-modelled years (i.e., 2051 to 2054) are not considered. 
    This is fine because the years based on which the annualised costs are calculated to align with the modelling years, so there are no more or fewer costs within the modelling horizon. 
    
    \item The total investment cost approach has a problem using the equation (\ref{eq:basic_total_cost_formulation}) in the objective function. It considers the total cost of the investment that covers the investment costs until 2054 but we only model until 2050.  
    To resolve this problem, we have to reformulate the total cost such that it only covers the modelled years. A standard approach is to consider the salvage value ($SV$).
\end{itemize}

The salvage value is what a company expects to get by selling or disassembling an asset at the end of its useful life, which will be the last modelled year in our case for these types of models. 
In the context of this document, the salvage value is the unit cost beyond the modelled years. Continue with our example, the unit cost from 2051 to 2054 is the difference between (\ref{eq:example_annualised_cost_formulation_ideal}) and (\ref{eq:example_annualised_cost_formulation_real}).

    \begin{equation}
     \label{eq:example_salvage_value}
     \begin{aligned}
                SV_{0} = & C_{0}^{A} \cdot (1 + \frac{1}{(1+WACC_0)^{11}} + \frac{1}{(1+WACC_0)^{12}} + \frac{1}{(1+WACC_0)^{13}} + \frac{1}{(1+WACC_0)^{14}}) 
     \end{aligned}
    \end{equation}

Generalizing (\ref{eq:example_salvage_value}) leads to (\ref{eq:salvage_value}).

\begin{equation}
 \label{eq:salvage_value}
    SV_{y} = C_{y}^{A} \cdot \sum_{j = Y+1}^{y+LT-1}\frac{1}{(1+WACC_{y})^{j - y}} \quad \forall y \in [0,Y]
\end{equation}

The salvage value $SV$ will depend on the year $y$, the last modelled year $Y$, the annualised value $C_{y}^{A}$, the lifetime $LT$, and the $WACC$. For instance, consider the following data for the first year ($y=0$):

\begin{itemize}
    \item Total capital cost $C_{0}^{T}=100$ [\euro/MW]
    \item The $WACC_{0}=5$ \%
    \item The last modelled year $Y=4$
    \item The economic lifetime $LT=8$ years
    \item Annualised capital cost using equation (\ref{eq:annualised_cost_from_total_cost}): $C_{0}^{A}=14.74$ [\euro/MW]
\end{itemize}

Therefore, using equation (\ref{eq:salvage_value}), the $SV_{0} = 33.01$ [\euro/MW] corresponds to the remaining unit cost in non-modelled years: 5, 6, and 7. The $SV$ is referenced to the first year $y=0$ and can be considered as a benefit of the investment in the objective function of the last modelled year which corresponds to the remaining cost in the non-modelled years, as in (\ref{eq:total_cost_formulation_with_salvage_value_a}). To put SV into perspective, if an investment is made in the last modelling year, after considering SV in the objection function, the total cost of this investment will only be the costs for this year, i.e., the same as the annualised cost. 

\begin{equation}
 \label{eq:total_cost_formulation_with_salvage_value_a}
    C^{I} = \sum_{y=0}^{Y-1}\frac{1}{(1+R)^{y}} \cdot (C_{y}^{T}-SV_{y})\cdot x_{y}
\end{equation}

\subsection{Milestone years} \label{sec:milestone_years_weight}
Salvage value can be used to cope with the modelling horizon ending before the lifetime ends. Now, we continue to discuss the other situation where milestone years are modelled instead of yearly investments.
When studying pathways to the year 2050, it is impractical to model yearly investment decisions on a large scale. A common approach is to use the so-called milestone years, where the major investment decisions will be made (e.g., years 2030, 2035, 2040, and 2050). 

Consider the information in Table \ref{tab:milestone_years_ex1} as an example. 
In the modelling horizon from 0 to 5, we only model the milestone years (i.e., years 0, 2, and 5) and we assume investments are made in these years. Years 1, 3, and 5 are not modelled, but their costs have to be incorporated. 

\begin{table}[ht!]
\centering
\begin{tabular}{|c|c|c|c|c|c|c|}
\hline
\textbf{Years $y$} & 0 & 1 & 2 & 3 & 4 & 5\\ 
\hline
\hline
\textbf{Milestone year $m$} & \textcolor{blue}{0} & - & \textcolor{blue}{2} & - & - & \textcolor{blue}{5}\\ 
\hline
\textbf{Milestone year weight $W_{m}^{I}$} & \textcolor{blue}{2} & - & \textcolor{blue}{3} & - & - & \textcolor{blue}{1} \\ 
\hline
\end{tabular}
\caption{Milestone years}
\label{tab:milestone_years_ex1}
\end{table}

The \textit{annualised investment cost} approach considers the costs of milestone years. Because the modelling horizon also includes non-modelled years, we need to incorporate those costs by using a weight coefficient $W_{m}^{I}$ (see equation (\ref{eq:annualised_investment_cost_milestone_weights})). By using this coefficient, we assume that the costs of the non-modelled years are represented by the costs of the modelled years. In this example, year 0 has a weight of 2 meaning that the cost of year 1 is the same as the cost of year 0 and the cost of year 0 covers the costs of two years (year 0 and year 1).

\begin{equation}
 \label{eq:annualised_investment_cost_milestone_weights}
  \begin{aligned}
    & \min_{x_{m},p_{mkt}} \quad C^{I} + C^{O}  \\
    & C^{I} = \sum_{m \in \mathcal{M}}\frac{1}{(1+R)^{m}} \cdot x_{m} \cdot \sum_{j = m | j \in \mathcal{M}}^{\min(m+LT-1,Y-1)}\frac{1}{(1+WACC_{m})^{j - m}} \cdot \textcolor{blue}{W_{j}^{I}} \cdot C_{m}^{A} \\ 
    & C^{O} = \sum_{m \in \mathcal{M}}\frac{1}{(1+R)^{m}} \cdot \textcolor{blue}{W_{m}^{I}} \cdot C_{m}^{op} \sum_{k}W_{mk}^{op} \cdot \sum_{t} p_{mkt}
  \end{aligned}
\end{equation} 

Please take note that the calculation for investment cost in equation (\ref{eq:annualised_investment_cost_milestone_weights}) only considers milestone years in the inner summation, owing to $j \in \mathcal{M}$.

The \textit{total investment cost} approach considers the total investment costs over both modelled and non-modelled years and only needs a weight coefficient in the operational costs $W_{y}^{I}$ to account for the non-modelled years. Moreover, to incorporate the investment costs beyond the modelled years, SV is used. Equation (\ref{eq:total_investment_cost_milestone_weights}) shows the formulation.

\begin{equation}
 \label{eq:total_investment_cost_milestone_weights}
  \begin{aligned}
    & \min_{x_{m},p_{mkt}} \quad C^{I} + C^{O}  \\
    & C^{I} = \sum_{m \in \mathcal{M}}\frac{1}{(1+R)^{m}} \cdot (C_{m}^{T}-SV_{m})\cdot x_{m} \\ 
    & C^{O} = \sum_{m \in \mathcal{M}}\frac{1}{(1+R)^{m}} \cdot \textcolor{blue}{W_{m}^{I}} \cdot C_{m}^{op} \sum_{k}W_{mk}^{op} \cdot \sum_{t} p_{mkt}
  \end{aligned}
\end{equation}    

Let us continue with the example in Table \ref{tab:milestone_years_ex1} (i.e., $\mathcal{M} = \{0,2,5\}$ and  $W_{m}^{I} = [2,3,1]$) and assume that $LT=6$. 
First, the investment costs using both methods are as follows.

\textit{Equation (\ref{eq:annualised_investment_cost_milestone_weights}) - annualised investment cost method:}
\begin{equation*}
  \begin{aligned}
     C^{I} = & \textcolor{orange}{C_{0}^{A}\left(2+\frac{3}{(1+WACC_{0})^{2}} +\frac{1}{(1+WACC_{0})^{5}}\right)}\cdot x_{0} \\ 
             & +\frac{1}{(1+R)^{2}} \cdot C_{2}^{A} \left(3+\frac{1}{(1+WACC_{2})^{3}} \right) \cdot x_{2} + \frac{1}{(1+R)^{5}}C_{5}^{A\cdot}x_{5} \\             
  \end{aligned}
\end{equation*}    

\textit{Equation (\ref{eq:total_investment_cost_milestone_weights}) - total investment cost method:}
\begin{equation*}
    C^{I} = \textcolor{orange}{C_{0}^{T}}x_{0} + \frac{1}{(1+R)^{2}} (C_{2}^{T}-SV_{2})\cdot x_{2} + \frac{1}{(1+R)^{5}} (C_{5}^{T}-SV_{5})\cdot x_{5}
\end{equation*}    

By simply inspecting the methods used, it is clear that the costs represented in the objective function are not the same as before when yearly investment decisions were made. For example, if we compare the first term in both costs (highlighted in \textcolor{orange}{orange}), we can see that the values differ according to equations (\ref{eq:total_cost_from_annualised_cost}) and (\ref{eq:annualised_cost_from_total_cost}).

Then, the operational cost is the same in both methods since the assumption here is that the milestone years represent the non-modelled years using the weight $W_{m}^{I}$. Applying the values of the example in equations (\ref{eq:annualised_investment_cost_milestone_weights}) and (\ref{eq:total_investment_cost_milestone_weights}), the operational cost is:

\begin{equation*}
  \begin{aligned}
    C^{O} = & \sum_{m \in \mathcal{M}}\frac{1}{(1+R)^{m}} \cdot \textcolor{blue}{W_{m}^{I}} \cdot C_{m}^{op} \sum_{k}W_{mk}^{op} \cdot \sum_{t} p_{mkt} \\
    & + 2 \cdot C_{0}^{op} \sum_{k}W_{0,k}^{op} \cdot \sum_{t} p_{0,kt} \\
    & + \frac{3}{(1+R)^{2}} \cdot C_{2}^{op} \sum_{k}W_{2,k}^{op} \cdot \sum_{t} p_{2,kt} \\
    & + \frac{1}{(1+R)^{5}} \cdot C_{5}^{op} \sum_{k}W_{5,k}^{op} \cdot \sum_{t} p_{5,kt} \\
  \end{aligned}
\end{equation*}

Accordingly, there are two issues with using the milestone method as shown in this subsection.

\begin{itemize}
    \item Because we assume that the costs of non-modelled years are the same with the modelled years, the time value of the non-modelled years is not considered. In this subsection, in the annualised approach, we assume that the investment costs of non-modelled years are the same with the modelled years before the non-modelled years, and thus, the investment costs are overestimated. For the operational costs, both methods have an overestimation.
    \item If the asset's lifetime ends between two milestone years, the $C^I$ and $C^O$ in the objective function cannot account for this and will assume the asset is available throughout the entire range between the milestone years. For instance, if $LT=4$ in Table \ref{tab:milestone_years_ex1}, the asset should be available until year 3. However, equations (\ref{eq:annualised_cost_from_total_cost}) will include them in the objective function until the final year (i.e., $y=5$), resulting in an additional two years in the objective function.
\end{itemize}

\section{Conclusions}

This document has outlined and compared the main approaches for modeling multi-year investments in energy system planning, with a particular focus on cost treatment and discounting mechanisms. Two commonly used representations of investment costs—the total cost (or overnight cost) approach and the annualised cost approach—have been presented both conceptually and mathematically. We demonstrated that while the two are mathematically equivalent under ideal conditions, their implementation may differ depending on data availability, modeling choices, and practical considerations.

We further analyzed two real-world modeling complications: (1) when the modeling horizon ends before the asset’s lifetime, and (2) when milestone years are used in place of annual investment decisions. To address the former, we introduced the concept of salvage value, which allows for an appropriate correction of total cost representations. To deal with the latter, we introduced milestone year weights, which approximate the influence of non-modelled years. The modeling logic and cost discounting mechanisms discussed here are implemented in the open-source \textit{TulipaEnergyModel} \cite{Tejada-Arango2023a}, and are documented in full detail in the accompanying technical reference \cite{Tejada-Arango2023b}.

While these adjustments make the models tractable for long-term studies, they also introduce simplifications that must be carefully considered. In particular, assumptions regarding the distribution and timing of investments and operations across milestone years can lead to over- or under-estimations of system costs. Therefore, modelers must strike a balance between computational efficiency and representation accuracy, depending on the purpose of the analysis and the availability of input data.

By providing a structured comparison and critical analysis of these methods, this document aims to support more transparent and informed choices in the design of multi-year investment models for energy system planning.

\bibliographystyle{unsrtnat}      
\bibliography{references}

\end{document}